%% file: order.tex
\documentclass{article}
\usepackage{a4wide}
\usepackage{theorem}
\usepackage{amssymb}
\input{define.tex}

\title{\vspace{-1.5cm}Ordering pure braid groups on closed surfaces}
\author{ Juan Gonz\'alez-Meneses}
\date{May, 2000}

\begin{document}
\maketitle

\begin{abstract}
We prove that the pure braid groups on closed, orientable surfaces are 
bi-orderable, and that the pure braid groups on closed, non-orientable
surfaces have generalized torsion, thus they are not bi-orderable.
\footnote{\hspace{-.65cm} Keywords: Braid - Surface -  
Orderable group.}
\footnote{\hspace{-.65cm} {\em Mathematics Subject Classification:} 
 Primary: 20F36. Secondary: 57N05.}
\footnote{\hspace{-.65cm} Partially supported by 
DGESIC-PB97-0723 and by the european network TMR Sing. Eq. Diff. et Feuill.}
\end{abstract}

\section{Introduction}

 The purpose of this paper is to answer the following question:  
Are pure braid groups on closed surfaces bi-orderable? We will prove
that the answer is positive for orientable surfaces, and negative
for the non-orientable ones. 

 In this section we give the basic definitions and 
classical results. We also explain what is known about orders on braid 
groups, and finally we state our results. In Section~2 we prove the 
stated result about closed, orientable surfaces. The 
non-orientable case is treated in Section~3.

\subsection{Braids on surfaces}

 Let $M$ be a compact surface, not necessarily orientable,
and let ${\cal P}=\{P_1,\ldots,P_n \}$ be a set of $n$ distinct 
points in $M$. Define a {\em $n$-braid based at ${\cal P}$} to be a 
collection $b = (b_1,\ldots,b_n)$ of disjoint smooth paths in 
$M\times [0,1]$, called {\em strings} of $b$, such that the $i$-th string 
$b_i$ runs monotonically in $t\in [0,1]$ from the point $(P_i,0)$ to some 
point $(P_j,1)$, $\; P_j\in {\cal P}$.

  An {\em isotopy} is defined as a deformation through braids,
which fixes the ends.  Multiplication of braids is defined by 
concatenation. The isotopy classes of braids with this multiplication 
form the group $B_n(M,{\cal P})$, called {\em braid group with $n$ 
strings on $M$ based at ${\cal P}$}. 
Note that the group $B_n(M,{\cal P})$ does not depend,
up to isomorphism, on the set ${\cal P}$ of points but only on the
cardinality $n=|{\cal P}|$. So we may write $B_n(M)$ in place of 
$B_n(M,{\cal P})$. Further details can be found in \cite{birman}.

 A {\em pure braid} is an element $b\in B_n(M)$ such that $b_i$ ends at 
$(P_i,1)$ for all $i=1,\ldots,n$. In other words, $b$ induces the trivial
permutation on ${\cal P}$. Pure braids form a normal subgroup of
$B_n(M)$ called  {\em pure braid group with $n$ 
strings on $M$}, and denoted by $PB_n(M)$. 

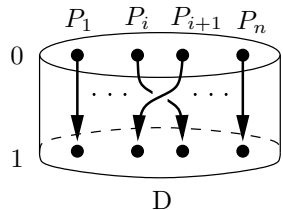
\begin{figure}[ht]
\centerline{\input{sigmai.pstex_t}}
\caption{The braid $\si_i\in B_n(D)$.}\label{sigmai}
\end{figure}

 Let $D$ be a closed disc. Notice that $B_n=B_n(D)$ is the classical braid 
group of Artin \cite{artin}. Recall that $B_n$
is generated by the set $\{\si_1,\ldots,\si_{n-1}\}$, where $\si_i$ is
the braid shown in Figure~\ref{sigmai}.

\vspace{.3cm}

\subsection{Orderable groups}

 A group $G$ is said to be {\em right-orderable} if there exists a strict
total ordering  $<$ on its elements which is invariant under 
right-multiplication: $g<h$ implies $gk < hk$ for all $g,h,k\in G$.
If $<$ is also invariant under left-multiplication, the $G$ is said to be 
{\em bi-orderable}.

 The following is a well known characterization of right-orderable 
and bi-orderable groups.

\begin{proposition}
 A group $G$ is right-orderable if and only if there exists a subset 
${\mathfrak p}\subset G$ such that ${\mathfrak p}^2 \subset {\mathfrak p}$ and 
$G= {\mathfrak p} \coprod \{ 1\} \coprod {\mathfrak p}^{-1}$. Moreover, $G$ is
bi-orderable if and only if there exists such ${\mathfrak p}$ also satisfying
$g {\mathfrak p} g^{-1}\subset {\mathfrak p}$, for all $g\in G$.
\end{proposition}

\begin{proof}
 If $G$ is right-orderable (bi-orderable), just take ${\mathfrak p}=
\{g\in G; \; 1<g\}$, the set of positive elements. Conversely, if there 
exists ${\mathfrak p}$ verifying the required hypothesis, then define the 
order $<$ by: $g<h$ if and only if $hg^{-1}\in {\mathfrak p}$.
\end{proof}

\vspace{.3cm}
 Let us state three properties of orderable groups, which can be found 
either in \cite{paris} or in \cite{rolfsenzhu}.

\begin{proposition}\label{proext}
 Let $1 \rightarrow A \stackrel{\al}{\rightarrow} B 
\stackrel{\be}{\rightarrow} C \rightarrow 1$ be an exact sequence of groups.
If $A$ and $C$ are right-orderable, with sets of positive elements 
${\mathfrak p}_A$ and ${\mathfrak p}_C$ respectively, then the set 
${\mathfrak p}_B= \al ({\mathfrak p}_A)\cup \be^{-1}({\mathfrak p}_C)$ defines a
right order on $B$. Moreover, if $A$ and $C$ are bi-orderable, the 
order on $B$ defined in such a way is a bi-order if and only if
$\al({\mathfrak p}_A)$ is normal in $B$ (that is, if the order on $A$ is 
preserved by conjugation in $B$). 
\end{proposition}

\begin{definition}
 We say that a group $G$ has {\em generalized torsion} if there exist
$g,h_1,\ldots,h_k\in G$, $g\neq 1$, such that 
$(h_1 g h_1^{-1})(h_2 g h_2^{-1})\cdots (h_k g h_k^{-1})=1$. 
\end{definition}

\begin{proposition}\label{progt}
A bi-orderable group has no generalized torsion (in particular, it is 
torsion-free).
\end{proposition}

%

\begin{proposition}
Let $G$ be a right-orderable group, and let $R$ be a ring with no zero 
divisors. Then the group ring $RG$ has no zero divisors. Moreover, the 
only units of $RG$ are the monomials $rg$, with $r$ invertible in $R$. 
\end{proposition}

\subsection{Ordering braid groups}

 Artin braid groups $B_n=B_n(D)$ are known to be right-orderable 
(\cite{dehornoy}, see also \cite{fenn}). However they are not 
bi-orderable if $n\geq 3$,
since they have generalized torsion \cite{neuwirth}. Nevertheless,
$P_n=PB_n(D)$ is bi-orderable, as shown in \cite{rolfsenzhu} and in
\cite{kimrolfsen}.

 In more generality, it is shown in \cite{rourkewiest} (see also 
\cite{shortwiest}) that the
mapping class groups of compact surfaces with boundary, fixing a finite
number of points, are right-orderable. Braid groups on compact surfaces
with boundary are subgroups of these mapping class groups, thus they are 
also right-orderable. 

 On the other hand, it is not known if the braid groups on closed surfaces 
are right-orderable or not. If $M$ is a closed surface different from
the sphere and from the projective plane, then $B_n(D)\subset B_n(M)$
(see \cite{parisrolfsen}). Hence, $B_n(M)$ also has generalized torsion,
and in conclusion, it cannot be bi-orderable. 

 In this paper, we will study the pure braid groups of closed surfaces.
We will show the following.

\begin{theorem}\label{teoori}
 If $M$ is a closed, orientable surface, then $PB_n(M)$ is bi-orderable.
\end{theorem}

\begin{theorem}\label{teonon}
 If $M$ is a closed, non-orientable surface, then $PB_n(M)$ has generalized
torsion, for $n\geq 2$. Therefore, $PB_n(M)$ is not bi-orderable
for $n\geq 2$.
\end{theorem}

\section{Closed, orientable surfaces}

 In this section we prove Theorem~\ref{teoori}. In Subsection~\ref{free}
we state that free groups and fundamental groups of orientable
surfaces are bi-orderable. In the case of free groups, we will explicitly 
define a bi-order. Then we see in Subsection~\ref{final} that 
$PB_n(M)$ is an extension of two groups, $K_n$ and $\pi_1(M)^n$, 
which are both bi-orderable. Moreover, the hypothesis of 
Proposition~\ref{proext} are satisfied, so $PB_n(M)$ turns to be
bi-orderable.

\subsection{Bi-order of free groups and fundamental groups}\label{free}

 We will explicitly define a bi-order on a given free group using the
so-called {\em Magnus expansion} \cite{magnus}. Let $F$ be a free group
with free system of generators $G=\{x_i \}_{i\in I}$ ($I$ not
necessarily finite). Let $\ZZ\lbra X_I \rbra$ be the ring of formal
power series over the non-commutative indeterminates $\{X_i\}_{i\in I}$.
The {\em Magnus expansion} of $F$ is a multiplicative homomorphism
$M:\: F\rightarrow \ZZ\lbra X_I \rbra$, such that
$ M(x_i)= 1+ X_i$ and  $ M(x_i^{-1})=1-X_i+X_i^2-\cdots$, 
for all $i\in I$. $M$ is known to be a well defined and
injective homomorphism, whose image is contained in 
$\{1+\eta \in \ZZ\lbra X_I \rbra;\; \eta(0)=0 \}$.

 Let us define a total order on $\ZZ\lbra X_I \rbra$. First, we choose
a total order on the set $\{X_i \}_{i\in I}$. Then, we order the monomials
of $\ZZ\lbra X_I \rbra$ as follows:
$$
  m_1 < m_2 \Leftrightarrow \left\{
\begin{array}{l}
 \mbox{deg}(m_1)< \mbox{deg}(m_2)  \\  
 \mbox{deg}(m_1)= \mbox{deg}(m_2)  \mbox{ and } m_1 <_{\mbox{lex}} m_2,
\end{array} \right.
$$
where deg means total degree (the sum of all exponents)  and  
$<_{\mbox{lex}}$ means smaller in the lexicographical order. 
The total order $\prec$ on $\ZZ\lbra X_I \rbra$ is defined as follows:
given $f,g\in \ZZ\lbra X_I \rbra$, we say that $f\prec g$ if and only if
the coefficient of the smallest non-trivial term of $g-f$ is positive. 

 Turning back to $F$, we define the {\em Magnus order} on it: given 
$a,b\in F$, we say that $a<b$ if and only if $M(a)\prec M(b)$. The
Magnus order is known to be a bi-order on $F$. Moreover,
one has:

\begin{theorem}{\em \cite{kimrolfsen}}\label{teoh1} 
The Magnus order on $F$ is preserved
under any $\Phi\in \mbox{\em Aut}(F)$ which induces the identity on
$H_1(F)=F/[F,F]$. $\fbox{}$
\end{theorem}

\vspace{.3cm} 
Let $\psi$ be a permutation of the set $\{X_i \}_{i\in I}$, and consider 
its extension $\Psi\in \mbox{\em Aut}(F)$. One has:

\begin{theorem}\label{teopsi}
 If $\psi$ preserves the order on $\{X_i \}_{i\in I}$, 
then $\Psi$ preserves the Magnus order on $F$.
\end{theorem}

\begin{proof}
 Notice that, under the action of such a $\Psi$, the degree and the 
lexicographical order on the monomials are preserved.
Hence, $\Psi$ preserves the order we defined on the monomials, thus the 
order $\prec$ on $\ZZ\lbra X_I \rbra$. Therefore, the Magnus order
on $F$ is also preserved.
\end{proof}

\vspace{.3cm}
 We finish this subsection with the following result.

\begin{theorem}{\em \cite{baumslag}}\label{teopi1} 
If $M$ is a closed, orientable surface, then $\pi_1(M)$ is a bi-orderable
group. $\fbox{}$
\end{theorem}

\subsection{$PB_n(M)$ is bi-orderable}\label{final}

 Let $M$ be a closed, orientable surface. Given a pure braid 
$b=(b_1,\ldots,b_n)\in PB_n(M)$, we can consider, for all $i=1,\ldots,n$,
the projection $\mu_i$ of $b_i$ over $M$. Since $b\in PB_n(M)$, $\mu_i$
is a loop in $M$ for all $i=1,\ldots,,n$, which represents an element 
of $\pi_1(M)$. This defines an epimorphism $\theta:\: PB_n(M)
\rightarrow \pi_1(M)^n$, which sends $(b_1,\ldots,b_n)$ to
$(\mu_1,\ldots,\mu_n)$ (see \cite{birman}).

 Define $K_n=\mbox{ker}(\theta)$. One has the exact sequence
$$
   1 \longrightarrow K_n \longrightarrow PB_n(M) 
\stackrel{\theta}{\longrightarrow} \pi_1(M)^n \longrightarrow 1.
$$

\begin{figure}[ht]
\centerline{\input{polairais.pstex_t}}
\caption{The polygon representing $M$ and the braids 
$a_{i,{2k+1}}$ and $a_{i,2k}$.}\label{polairais}
\end{figure}
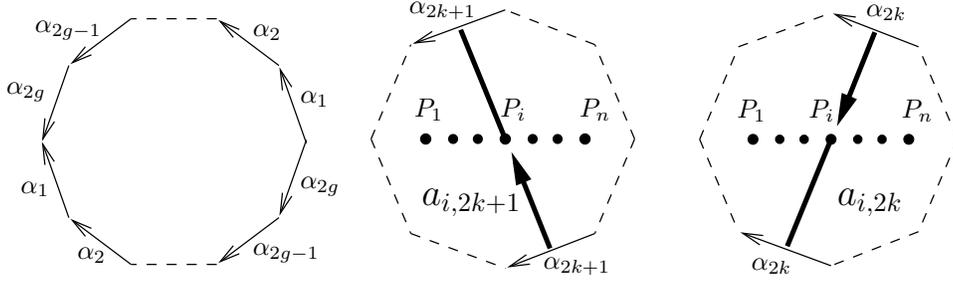

Let us prove that $K_n$  is bi-orderable. First, we represent $M$ as
a polygon of $4g$ sides, identified in the way of Figure~\ref{polairais}.
For all $i=1,\ldots,n$ and all $r=1,\ldots,2g$, we define the braid 
$a_{i,r}\in PB_n(M)$ as in Figure~\ref{polairais}:
The $i$-th string of $a_{i,r}$ is 
$(s_{i,r}(t),t)\in M \times [0,1]$,
where $s_{i,r}$ is a loop in $M$ based at $P_i$ which goes through the wall
$\al_r$; it goes upwards if $r$ is odd and downwards if $r$ is even. 
The $j$-th string of $a_{i,r}$ is $(P_j,t)$ (the trivial string) for all 
$j\neq i$.

 Let now $\Om=\{\om_1,\ldots,\om_{2g} \}$ be a set of generators of 
$\pi_1(M)$, where $g$ is the genus of $M$. Take $\Om$ in such a way 
that $\pi_1(M)=\left< \Om; \: \om_1 \cdots \om_{2g} \om_1^{-1}\cdots 
\om_{2g}^{-1}=1 \right>$. For all $\ga\in \pi_1(M)$, choose
a unique word $\widetilde{\ga}$ over $\Om \cup \Om^{-1}$ representing
$\ga$. We denote by $\widetilde{\ga}_{(i)}$ the pure braid obtained 
from $\widetilde{\ga}$ by replacing $\om_r^{\pm 1}$ with $a_{i,r}^{\pm 1}$.
Now, for all $i,j\in \{1,\ldots,n \}$, $i<j$, we define the braid
$$
  t_{i,j}=t_{j,i}= \si_i \cdots \si_{j-2} \si_{j-1}^2 \si_{j-2}^{-1}
 \cdots  \si_i^{-1} \in PB_n(M).
$$
Finally, for all $i,j\in \{1,\ldots,n \}$, $i\neq j$, and all 
$\ga\in \pi_1(M)$, we define $f_{i,j,\ga}= \widetilde{\ga}_{(i)}
 t_{i,j} \widetilde{\ga}_{(i)}^{-1}$.

\begin{theorem}{\em \cite{gmparis}}\label{teodecomp}
 One has $K_n=( F_n \rtimes(F_{n-1}\rtimes(\cdots(F_3\rtimes F_2)\cdots)))$,
where for all $i=1,\ldots,n-1$, $F_{(n+1)-i}$ is the free group 
freely generated by ${\cal F}_{i,n}=\{f_{i,j,\ga};\; i< j\leq n, \:
 \ga\in \pi_1(M)\}$. Moreover, for all $m=2,\ldots,n-1$,
$K_m=(F_m \rtimes(\cdots(F_3\rtimes F_2)\cdots ))$ acts trivially on
$H_1(F_{m+1})$. $\fbox{}$ 
\end{theorem}

\begin{corollary}\label{cor}
 $K_n$ is bi-orderable.
\end{corollary}

\begin{proof}
 We argue by induction on $n$. If $n=2$, then $K_n=F_2$ is a free group
(of infinite rank), so it is bi-orderable. Suppose then that $n>2$,
and that $K_{n-1}$ is bi-orderable. By Theorem~\ref{teodecomp}, we 
have an exact sequence
$$
  1 \longrightarrow F_n \longrightarrow K_n \longrightarrow K_{n-1}
\longrightarrow 1,
$$
where $K_n= F_n \rtimes K_{n-1}$. By definition of bi-order, 
conjugation by an element of $F_n$ is an automorphism of $F_n$ which 
preserves the Magnus order. We also know, by Theorem~\ref{teodecomp},
that conjugation by an element of $K_{n-1}$ is 
an automorphism of $F_n$ which is trivial on $H_1(F_n)$. 
Hence, by Theorem~\ref{teoh1}, it also preserves the Magnus 
order on $F_n$. Therefore, conjugation by an element of $K_n$ preserves
the Magnus order of $F_n$ and thus, by Proposition~\ref{proext}, $K_n$ is 
bi-orderable.
\end{proof}

\vspace{.3cm}
 Let us define an explicit bi-order on $K_n$. First, for all 
$i=1,\ldots,n-1$, we order ${\cal F}_{i,n}$ as follows:
$$
  f_{i,j,\ga}< f_{i,k,\de} \Leftrightarrow \left\{
\begin{array}{l}
 j<k  \\   j=k \mbox{ and } \ga <_{\pi_1} \de,
\end{array} \right.
$$ 
where $<_{\pi_1}$ is a fixed bi-order of $\pi_1(M)$. Then, we consider 
the Magnus order on each $F_{(n+1)-i}$ corresponding to this order on 
${\cal F}_{i,n}$. The bi-order on $K_n$ which yields from 
Corollary~\ref{cor} is the following: for $k,k'\in K_n$, write
$k=k_1k_2\cdots k_{n-1}$ and $k'=k'_1k'_2\cdots k'_{n-1}$,
where $k_i,k'_i\in F_{(n+1)-i}$. Then $k<k'$ if and only
if $k_j<k'_j$ for the greatest $j$ such that $k_j\neq k'_j$.

\vspace{.3cm}
\noindent {\sc Proof of Theorem~\ref{teoori}}: The direct product of 
bi-orderable groups is clearly bi-orderable, hence, by
Theorem~\ref{teopi1}, $\pi_1(M)^n$ is bi-orderable. So, by
Proposition~\ref{proext}, we only need to show that conjugation by
an element of $PB_n(M)$ is an automorphism of $K_n$ which preserves
the order. 

 Conjugation by an element of $K_n$ preserves the order by definition
of bi-order. Hence, it suffices to show the above claim for the conjugation
by pre-images under $\theta$ of the generators of $\pi_1(M)^n$. A set of
such pre-images is $\{a_{i,r}; \; i=1,\ldots,n,\; r=1,\ldots,2g \}$.
Now, in [G-MP, Lemma~3.15] it is shown that the following relations hold in
$H_1(K_n)$:
$$
    a_{i,r}f_{j,k,\ga}a_{i,r}^{-1}\equiv \left\{ 
\begin{array}{ll}
 f_{j,k,\ga}  & \mbox{if } i\neq j,k \\
 f_{j,k,(\om_r\ga)} & \mbox{if }  i=j \\
f_{j,k,(\ga \om_r^{-1})} \hspace{1cm} & \mbox{if }  i=k. 
\end{array} \right.
$$
We claim that the action of $a_{i,r}$ preserves the Magnus order
on each $F_m$, $m=2,\ldots,n-1$, and hence, it preserves the order on 
$K_n$. Clearly, the action of $a_{i,r}$ on $K_n$ is the composition of
an automorphism $\Psi_{i,r}$ which permutes the generators of each $F_m$,
with an automorphism $\Phi_{i,r}$ which is trivial on $H_1(K_n)$.
Therefore, by Theorems \ref{teoh1} and \ref{teopsi}, it suffices
to prove that the permutation induced by $\Psi_{i,r}$ on 
${\cal F}_{j,n}$ $(j=1,\ldots,n-1)$ preserves the defined order
on ${\cal F}_{j,n}$.

 Let then $f_{j,k,\ga},f_{j,l,\de}\in {\cal F}_{j,n}$, where
$f_{j,k,\ga} < f_{j,l,\de}$.

\noindent {\bf Case 1:} If $k<l$, then
$\Psi_{i,r}(f_{j,k,\ga})=f_{j,k,\ga'}<f_{j,l,\de'}=\Psi_{i,r}(f_{j,l,\de})$,
where $\ga'$ and $\de'$ are determined by the above relations.

\noindent {\bf Case 2:} If $k=l$ and $\ga <_{\pi_1} \de$, then there are
three possibilities. First, if $i\neq j,k$, one has 
$\Psi_{i,r}(f_{j,k,\ga})=f_{j,k,\ga}<f_{j,k,\de}=\Psi_{i,r}(f_{j,k,\de})$.
If $i=j$, one has $\Psi_{i,r}(f_{j,k,\ga})=f_{j,k,(\om_r\ga)}<
f_{j,k,(\om_r\de)} =\Psi_{i,r}(f_{j,k,\de})$, since 
$\om_r\ga <_{\pi_1} \om_r\de$ ($<_{\pi_1}$ is a left-order). 
Finally, if $i=k$, then
$\Psi_{i,r}(f_{j,k,\ga})=f_{j,k,(\ga\om_r^{-1})}<f_{j,k,(\de\om_r^{-1})}=
\Psi_{i,r}(f_{j,k,\de})$, since $\ga\om_r^{-1} <_{\pi_1} \de\om_r^{-1}$ 
($<_{\pi_1}$ is a right-order).

 Therefore, $PB_n(M)$ is a bi-orderable group. $\fbox{}$

\section{Closed, non-orientable surfaces}

 We turn now to prove Theorem~\ref{teonon}. Let $M$ be a closed, 
non-orientable surface, and let ${\cal P}=\{P_1,\ldots, P_n \}\subset M$.
Then, there exists a submanifold $N$, $\: {\cal P}\subset N \subset M$,
such that $N$ is homeomorphic to a M\"obius strip. 

 Consider the subset $C= \RR \times [0,1] \times [0,1] \subset \RR^3$,
and identify $(x,y,t)\sim (x+1,1-y,t)$ for all $x\in \RR$, and all 
$y,t\in [0,1]$. One has $(C/\sim) \simeq N\times [0,1]$.
We choose $P_i=(1/2,p_i)\in N$, where $p_i=\frac{i}{n+1}$ for all 
$i=1,\ldots,n$. 

\begin{figure}[ht]\label{gade}
\centerline{\input{gade.pstex_t}}
\caption{The braids $\Ga, \De \in B_n(N)$.}
\end{figure}
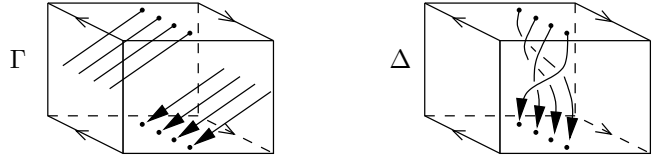

 Denote by $\Ga=(\ga_1,\ldots,\ga_n)$ the braid on $N$
defined as follows: $\ga_i(t)=(1/2-t,p_i,t)$, for all $i=1,\ldots,n$.
It is the braid represented in Figure~\ref{gade}, for $n=4$. 
Notice that $\Ga$ is not a pure braid (going through the wall reverses 
the orientation). Denote by $\De$ the following braid on $N$: 
$\De=(\si_1\cdots \si_{n-1})(\si_1\cdots \si_{n-2})\cdots(\si_1\si_2)\si_1$.
It is also drawn in Figure~\ref{gade}, for $n=4$.
Remark that we represent $N\times\{ 0\}$ above
$N\times\{ 1\}$ to agree with the usual orientation of braids 
(pointing downwards).

\vspace{.3cm}

 Now, given a point $p=(x,y,t)\in C/\sim$, we denote by $\overline{p}$
the image of $p$ under the symmetry of plane $y=1/2$. That is,
$\overline{p}=(x,1-y,t)$. In the same way, given a braid $b\in B_n(N)$,
we denote by $\overline{b}$ the image of $b$ under the same symmetry
($\overline{b}_i(t)=\overline{b_i(t)}$). 

\begin{lemma}
 For all $b\in B_n(N)$, one has $\Ga b \Ga^{-1}= \overline{b}$.
\end{lemma}

\begin{proof}
 Let $b=(b_1,\ldots,b_n)\in B_n(N)$. For all $i=1,\ldots,n$, we write
$b_i(t)=(\be_i(t),t)$, where $\be_i$ is a path on $N$. We denote by
$\ep$ the permutation induced by $b$ on ${\cal P}$. Consider the 
braid  $c=(c_1,\ldots,c_n)=\Ga b \Ga^{-1}$. For all $i=1,\ldots,n$, one has
$$
  c_i(t)= \left\{
\begin{array}{ll}
  (1/2-3t, p_i, t)  &  \mbox{if } t\in [0, 1/3],  \\
  (\be_{n+1-i}(3t-1),t)  &  \mbox{if } t\in [1/3, 2/3],  \\
  (1/2+(3t-2), p_{\ep(n+1-i)}, t) \hspace{1cm} &  \mbox{if } t\in [2/3, 1].  
\end{array} \right.
$$
Now consider the isotopy 
$H:\: (N\times[0,1])\times [0,1] \rightarrow N\times[0,1]$ given by
$$
    H(x,y,t,s)= \left\{
\begin{array}{ll}
  (x+3ts, y, t)  &  \mbox{if } t\in [0, 1/3],  \\
  (x+s, y, t)  &  \mbox{if } t\in [1/3, 2/3],  \\
  (x+3(1-t)s, y, t) \hspace{1cm} &  \mbox{if } t\in [2/3, 1],  
\end{array} \right.
$$
Then $H(c,0)=c$ and $H(c,1)=1\: \overline{b} \: 1 \simeq \overline{b}$.
\end{proof}

\vspace{.3cm}
\noindent {\sc Proof of Theorem~\ref{teonon}}: Take a closed disc $D$,
$\; {\cal P}\subset D\subset N\subset M$. It is well known that, in
$B_n(D)$, $\; \De \si_i \De^{-1}= \si_{n-i}$ for all $i=1,\ldots,n-1$.
So, the same relation holds in $B_n(N)$ and in $B_n(M)$ (every isotopy can
be extended by the identity outside $D$). Hence, for $i=1,\ldots,n-1$,
one has $(\Ga \De)\si_i (\Ga \De)^{-1}=\Ga \: \si_{n-i}\: \Ga^{-1}=
\overline{\si_{n-i}}=\si_i^{-1}$ in $B_n(N)$, thus in $B_n(M)$. 

 Therefore, one has $\si_i \left[(\Ga \De)\si_i (\Ga \De)^{-1} \right]= 1$,
and so, $\si_i^2 \left[(\Ga \De)\si_i^2 (\Ga \De)^{-1} \right]= 1$
in $B_n(M)$. It suffices to notice that $\si_i^2$ and $\Ga\De$ are pure
braids, and that $\si_i^2\neq 1$ in $PB_n(M)$ for all $i=1,\ldots,n$, 
where $n\geq 2$, to conclude that $PB_n(M)$ has generalized torsion.
$\fbox{}$

%

\begin{tabular}{ll}
\noindent J. GONZ\'ALEZ-MENESES \\
Universit\'e de Bourgogne \hspace{4truecm} & Departamento de \'Algebra \\
Laboratoire de Topologie & Facultad de Matem\'aticas\\
UMR 5584 du CNRS & Universidad de Sevilla\\
B. P. 47870  & C/ Tarfia, s/n\\
21078 - Dijon Cedex (France)&  41012 - Sevilla (Spain)\\
{\em jmeneses@u-bourgogne.fr}& {\em meneses@algebra.us.es}
\end{tabular}

\end{document}

%% file: define.tex
\newcommand\al{\alpha}
\newcommand\be{\beta}
\newcommand\ga{\gamma}
\newcommand\de{\delta}
\newcommand\ep{\varepsilon}

\newcommand\si{\sigma}

\newcommand\om{\omega}

\newcommand\Ga{\Gamma}
\newcommand\De{\Delta}

\newcommand\Om{\Omega}

\newcommand\lbra{[\hspace{-2pt}[}
\newcommand\rbra{]\hspace{-2pt}]}

\def\RR{\mathbb R}
\def\ZZ{\mathbb Z}

  \newcounter{numero}[section]
  \newcounter{aux1}
  \newcounter{aux2}
  \renewcommand{\thenumero}{\thesection.\arabic{numero}}
  \newenvironment{corollary}{
  \dimen255=\parindent \parindent=0in
  \refstepcounter{numero}{\vspace{.3cm}\bf Corollary \thenumero .}\
  \it}{\rm \parindent=\dimen255\par}

  \newenvironment{definition}{
  \dimen255=\parindent \parindent=0in
  \refstepcounter{numero}{\vspace{.3cm}\bf Definition \thenumero .}\
  }{\rm \parindent=\dimen255\par}

  \newenvironment{theorem}{
  \dimen255=\parindent \parindent=0in
  \refstepcounter{numero}{\vspace{.3cm}\bf Theorem \thenumero .}\
  \it}{\rm \parindent=\dimen255\par}

  \newenvironment{lemma}{
  \dimen255=\parindent \parindent=0in
  \refstepcounter{numero} {\vspace{.3cm}\bf Lemma \thenumero .}\
  \it}{\rm \parindent=\dimen255\par}

  \newenvironment{proposition}{
  \dimen255=\parindent \parindent=0in
  \refstepcounter{numero} {\vspace{.3cm}\bf Proposition \thenumero .}\ \it
  }{\rm \parindent=\dimen255\par}

  \newenvironment{proof}{
  {\vspace{.3cm}\noindent\sc{Proof:}}\
  }{$\fbox{}$ \par}

%% file: sigmai.pstex_t
\begin{picture}(0,0)%
\special{psfile=sigmai.pstex}%
\end{picture}%
\setlength{\unitlength}{3947sp}%
\begingroup\makeatletter\ifx\SetFigFont\undefined
\def\x#1#2#3#4#5#6#7\relax{\def\x{#1#2#3#4#5#6}}%
\expandafter\x\fmtname xxxxxx\relax \def\y{splain}%
\ifx\x\y   
\gdef\SetFigFont#1#2#3{%
  \ifnum #1<17\tiny\else \ifnum #1<20\small\else
  \ifnum #1<24\normalsize\else \ifnum #1<29\large\else
  \ifnum #1<34\Large\else \ifnum #1<41\LARGE\else
     \huge\fi\fi\fi\fi\fi\fi
  \csname #3\endcsname}%
\else
\gdef\SetFigFont#1#2#3{\begingroup
  \count@#1\relax \ifnum 25<\count@\count@25\fi
  \def\x{\endgroup\@setsize\SetFigFont{#2pt}}%
  \expandafter\x
    \csname \romannumeral\the\count@ pt\expandafter\endcsname
    \csname @\romannumeral\the\count@ pt\endcsname
  \csname #3\endcsname}%
\fi
\fi\endgroup
\begin{picture}(1621,1325)(911,-656)
\put(911,-391){\makebox(0,0)[lb]{\smash{\SetFigFont{10}{12.0}{rm}1}}}
\put(1801,-656){\makebox(0,0)[lb]{\smash{\SetFigFont{10}{12.0}{rm}D}}}
\put(1251,489){\makebox(0,0)[lb]{\smash{\SetFigFont{10}{12.0}{rm}$P_1$}}}
\put(1626,504){\makebox(0,0)[lb]{\smash{\SetFigFont{10}{12.0}{rm}$P_i$}}}
\put(1916,504){\makebox(0,0)[lb]{\smash{\SetFigFont{10}{12.0}{rm}$P_{i+1}$}}}
\put(2331,479){\makebox(0,0)[lb]{\smash{\SetFigFont{10}{12.0}{rm}$P_n$}}}
\put(911,259){\makebox(0,0)[lb]{\smash{\SetFigFont{10}{12.0}{rm}0}}}
\end{picture}

%% file: polairais.pstex_t
\begin{picture}(0,0)%
\special{psfile=polairais.pstex}%
\end{picture}%
\setlength{\unitlength}{3947sp}%
\begingroup\makeatletter\ifx\SetFigFont\undefined
\def\x#1#2#3#4#5#6#7\relax{\def\x{#1#2#3#4#5#6}}%
\expandafter\x\fmtname xxxxxx\relax \def\y{splain}%
\ifx\x\y   
\gdef\SetFigFont#1#2#3{%
  \ifnum #1<17\tiny\else \ifnum #1<20\small\else
  \ifnum #1<24\normalsize\else \ifnum #1<29\large\else
  \ifnum #1<34\Large\else \ifnum #1<41\LARGE\else
     \huge\fi\fi\fi\fi\fi\fi
  \csname #3\endcsname}%
\else
\gdef\SetFigFont#1#2#3{\begingroup
  \count@#1\relax \ifnum 25<\count@\count@25\fi
  \def\x{\endgroup\@setsize\SetFigFont{#2pt}}%
  \expandafter\x
    \csname \romannumeral\the\count@ pt\expandafter\endcsname
    \csname @\romannumeral\the\count@ pt\endcsname
  \csname #3\endcsname}%
\fi
\fi\endgroup
\begin{picture}(5695,1784)(406,-1732)
\put(5446,-773){\makebox(0,0)[lb]{\smash{\SetFigFont{10}{12.0}{rm}$P_{i}$}}}
\put(1955,-262){\makebox(0,0)[lb]{\smash{\SetFigFont{10}{12.0}{rm}$\al_2$}}}
\put(853,-1661){\makebox(0,0)[lb]{\smash{\SetFigFont{10}{12.0}{rm}$\al_2$}}}
\put(1955,-1631){\makebox(0,0)[lb]{\smash{\SetFigFont{10}{12.0}{rm}$\al_{2g-1}$}}}
\put(2252,-1214){\makebox(0,0)[lb]{\smash{\SetFigFont{10}{12.0}{rm}$\al_{2g}$}}}
\put(2252,-679){\makebox(0,0)[lb]{\smash{\SetFigFont{10}{12.0}{rm}$\al_1$}}}
\put(2937,-113){\makebox(0,0)[lb]{\smash{\SetFigFont{10}{12.0}{rm}$\al_{2k+1}$}}}
\put(3781,-1720){\makebox(0,0)[lb]{\smash{\SetFigFont{10}{12.0}{rm}$\al_{2k+1}$}}}
\put(586,-232){\makebox(0,0)[lb]{\smash{\SetFigFont{10}{12.0}{rm}$\al_{2g-1}$}}}
\put(406,-619){\makebox(0,0)[lb]{\smash{\SetFigFont{10}{12.0}{rm}$\al_{2g}$}}}
\put(496,-1244){\makebox(0,0)[lb]{\smash{\SetFigFont{10}{12.0}{rm}$\al_1$}}}
\put(3511,-768){\makebox(0,0)[lb]{\smash{\SetFigFont{10}{12.0}{rm}$P_{i}$}}}
\put(4006,-768){\makebox(0,0)[lb]{\smash{\SetFigFont{10}{12.0}{rm}$P_{n}$}}}
\put(2971,-768){\makebox(0,0)[lb]{\smash{\SetFigFont{10}{12.0}{rm}$P_{1}$}}}
\put(3016,-1321){\makebox(0,0)[lb]{\smash{\SetFigFont{14}{16.8}{rm}$a_{i,2k+1}$}}}
\put(4996,-773){\makebox(0,0)[lb]{\smash{\SetFigFont{10}{12.0}{rm}$P_{1}$}}}
\put(6031,-773){\makebox(0,0)[lb]{\smash{\SetFigFont{10}{12.0}{rm}$P_{n}$}}}
\put(5806,-133){\makebox(0,0)[lb]{\smash{\SetFigFont{10}{12.0}{rm}$\al_{2k}$}}}
\put(5086,-1726){\makebox(0,0)[lb]{\smash{\SetFigFont{10}{12.0}{rm}$\al_{2k}$}}}
\put(5626,-1321){\makebox(0,0)[lb]{\smash{\SetFigFont{14}{16.8}{rm}$a_{i,2k}$}}}
\end{picture}

%% file: gade.pstex_t
\begin{picture}(0,0)%
\special{psfile=gade.pstex}%
\end{picture}%
\setlength{\unitlength}{3947sp}%
\begingroup\makeatletter\ifx\SetFigFont\undefined
\def\x#1#2#3#4#5#6#7\relax{\def\x{#1#2#3#4#5#6}}%
\expandafter\x\fmtname xxxxxx\relax \def\y{splain}%
\ifx\x\y   
\gdef\SetFigFont#1#2#3{%
  \ifnum #1<17\tiny\else \ifnum #1<20\small\else
  \ifnum #1<24\normalsize\else \ifnum #1<29\large\else
  \ifnum #1<34\Large\else \ifnum #1<41\LARGE\else
     \huge\fi\fi\fi\fi\fi\fi
  \csname #3\endcsname}%
\else
\gdef\SetFigFont#1#2#3{\begingroup
  \count@#1\relax \ifnum 25<\count@\count@25\fi
  \def\x{\endgroup\@setsize\SetFigFont{#2pt}}%
  \expandafter\x
    \csname \romannumeral\the\count@ pt\expandafter\endcsname
    \csname @\romannumeral\the\count@ pt\endcsname
  \csname #3\endcsname}%
\fi
\fi\endgroup
\begin{picture}(3842,924)(226,-523)
\put(2611, 29){\makebox(0,0)[lb]{\smash{\SetFigFont{10}{12.0}{rm}$\De$}}}
\put(226, 29){\makebox(0,0)[lb]{\smash{\SetFigFont{10}{12.0}{rm}$\Ga$}}}
\end{picture}